\documentstyle[12pt,a4]{article}

\newcommand{\be}{\begin{equation}}
\newcommand{\ee}{\end{equation}}
\newcommand{\bea}{\begin{eqnarray}}
\newcommand{\eea}{\end{eqnarray}}

\title{The hidden symmetry algebras of a class of quasi-exactly 
solvable multi dimensional operators}

\author
{Y. Brihaye{\thanks{E-mail : Yves.Brihaye@umh.ac.be}} 
\ and J. Nuyts{\thanks{E-mail : Jean.Nuyts@umh.ac.be.. Work supported in
part by the Belgian Fonds National de la Recherche Scientifique.}}
 \\
Department of Mathematical Physics\\
University of Mons\\
Av. Maistriau, B-7000 MONS, Belgium}
\begin{document}

\begin{titlepage}
\maketitle
\thispagestyle{empty}
\begin{abstract}
Let $P(N,V)$ denote the vector space of polynomials 
of maximal degree less than or equal to $N$ in $V$
independent variables. This space
is preserved by the enveloping algebra generated by a set
of linear, differential operators representing
the Lie algebra $gl(V+1)$. We establish the 
counterpart of this property for the vector space
$P(M,V) \oplus P(N,V)$ for any values of the integers $M,N,V$. 
We show that the operators preserving $P(M,V) \oplus P(N,V)$ generate
an abstract superalgebra (non linear if $\Delta=\mid M-N\mid\geq 2$).
A family of algebras is also constructed, extending this particular
algebra by $\Delta -1$ arbitrary complex parameters. 
\end{abstract}
\vfill
\end{titlepage}

\section{Introduction}
Quasi-exactly solvable (QES) equations refer to a class of spectral  
differential equations for which a part of the 
spectrum can be obtained by solving 
algebraic equations \cite{tur0,ush,shif}.
Linear differential operators preserving a finite dimensional
space of smooth functions constitute in this respect a basic 
ingredient in the topic of QES equations.    

In the case of operators of one real variable  acting on 
a scalar function, the possibilities of finite dimensional
invariant vector spaces are rather limited \cite{tur1}.
Up to a change of the variable and a redefinition of the function,
the vector space can only be the set $P(N)$ of polynomials with degree less
than or equal to $N$. 
The relevant operators are the elements of
the enveloping algebra of $sl(2)$ whose generators are 
suitably represented by three differential operators \cite{tur0,tur1}.
The QES equations they define therefore possess an $sl(2)$ hidden symmetry.

When the number of variables or (and) the number of components
of the function is (are) larger than one,
the number of possible invariant vector spaces 
and  hidden symmetry algebras increases considerably.
The  scalar QES operators in two variables were 
classified in \cite{gko1,gko2}. Seven inequivalent 
spaces of functions appear to be possible. Correspondingly, 
the hidden symmetry algebras can be of several
types, e.g. $sl(2)$, $sl(2) \otimes sl(2)$, $sl(3)$. 
A few cases of this classification
generalise easily to operators involving 
an arbitrary number of variables. 
In particular, the case labelled 2.3 in ref. \cite{gko2}
can be extended  to the space  $P(N,V)$  
of polynomials of maximal degree less than or equal to $N$
in their $V$ independent variables.
The related algebra is  $sl(V+1)$.

The construction of the matrix operators in one variable
preserving the direct sum $P(N_1) \oplus \ldots \oplus P(N_k)$ 
has also been considered \cite{bk,bggk}. These operators
are closely related to graded algebras. As an example,
the case $P(N) \oplus P(N+1)$ is related to the graded Lie
algebra $osp(2,2)$ \cite{shtu}. 

The purpose of this paper 
is to classify  the operators preserving the vector space 
$P(M,V)$ \break 
$\oplus P(N,V)$ for arbitrary values of the integers $M,N,V$  
and to construct a series of associative algebras 
corresponding to the hidden symmetries of these operators.
In Sect. 2 we fix the notations and point out the relevant  
representations of the algebra $gl(n)$.
The $2 \times 2$ matrix operators preserving the space  
$P(N,V) \oplus P(M,V)$ are contructed in Sect. 3 and are
shown to obey a set of normal ordering rules. 
In Sect. 4, these ordering rules are modified into sets 
of commutation and anticommutation relations which fulfil
all Jacobi identities. We obtain in this way a series 
of  associative abstract algebras which appear to be labelled by 
$V$, by  $\Delta=\mid N-M\mid $ and by $\Delta -1$ arbitrary complex
parameters.   
The technical details related to the proof of our main result
are given in Sect. 5.

\section{Operators preserving $P(N,V)$}

\par Let $N,V$ be two positive integers. Let $x_i\ (i=1,\ldots,V)$ 
represent $V$ independent real variables. We define the finite dimensional
vector space $P(N,V)$ of polynomials in the variables $x_i$ and
of maximal total degree  $N$
\be
P(N,V) = {\rm span} \ \lbrace x^{n_1}_1, x^{n_2}_2\ldots x^{n_V}_V
\rbrace
\quad , \quad  0 \leq \sum_{j=1}^V n_j \leq N 
\ee
\be
P(N,1) \equiv P(N) = {\rm span} \ \lbrace 1,x,\ldots,x^N \rbrace \ .
\ee
The dimension of $P(N,V)$ is given by
\be
 1 + V + {V(V+1)\over 2} + \ldots + 
{V(V+1) \ldots (V+N-1) \over {N!}} = C^{N+V}_V \ .
\ee

The set of linear differential operators preserving $P(N,V)$ can be perceived
as the enveloping algebra generated by the following operators 
\begin{eqnarray}
J^0_0(N) &=& D-N \quad , \quad 
D\equiv \sum^V_{j=1} x_j {\partial \over {\partial x_j}} \nonumber \\
J^k_0(N) &=& {\partial \over
{\partial x_k}}\quad , \quad k=1,\ldots, V   \nonumber \\
J^0_k(N) &=& -x_k(D-N)\quad , \quad k=1,\ldots, V \nonumber \\
J^l_k(N) &=& - x_k{\partial\over{\partial x_l}}\quad , \quad k,l=1,\ldots, V
\ .
\label{j}
\end{eqnarray}
These $(V+1)^2$ independent operators fulfil the commutation 
rules of the Lie algebra $gl(V+1)$. 
Acting on the finite dimensional space $P(N,V)$, they
lead to an irreducible representation of this algebra.
The commutation relations are 
\be
[J^{b}_{a}, J^{d}_{c}] = 
 \delta^{d}_{a} J^{b}_{c}
- \delta^{b}_{c}J^{d}_{a} 
\quad \ \ , \ \ a,b,c,d = 0,1, \ldots , V \ .
\label{comj}
\ee

Within the representation (\ref{j}), the Casimir operators of $gl(V+1)$   
\be 
C_p \equiv \sum_{a_1,\ldots, a_p = 0}^V 
 J^{a_1}_{a_2}  J^{a_2}_{a_3} \ldots  J^{a_p}_{a_1}
\ \ \ , \ \ \
 p=1, \ldots, V+1
\label{cas}
\ee
have the values   
$C_p= (-1)^p N(N+V)^{p-1}$.
The operators ${J'}_a^b$ defined by 
\be
   {J'}_{a}^{b} =   J_{a}^{b} 
                    + C \delta_{a}^{b}  \  \ ,
\label{trans}
\ee   
where $C$ is any operator which commutes with all $J$'s, 
satisfy also the relations (\ref{comj}).
For instance, this is the case for the $(V+1)^2-1$ independent operators 
\be
\tilde J_{a}^{b}=J_{a}^{b}- {1 \over V+1} C_1 \delta_{a}^{b}
\ee
(since $\tilde C_1 \equiv 0$) which 
form  an irreducible 
representation of $sl(V+1)$ when acting on $P(N,V)$.
The usual form \cite{tur0} of the of $sl(2)$  generators    
\begin{eqnarray}
J_+(N) &=& -\tilde J_1^0 = x(x\partial _x-N) \nonumber \\
J_0(N) &=& - \tilde J^1_1 = (x\partial _x -{N\over 2}) \nonumber \\
J_-(N) &=& \tilde J_0^1 = \partial _x
\end{eqnarray}
is recovered for  $V=1$.
These operators play a major role in the topic of 
quasi-exactly solvable equations.

\par More generally, an element, say $A$,
of the enveloping algebra constructed over the 
$J^{b}_{a}$ (or the $\tilde J^{b}_{a}$) is a
quasi-exactly solvable operator preserving $P(N,V)$.
That is to say that  the spectral equation
\be
A p = \lambda p \quad , \quad p \in P(N,V)
\ee
admits $C^{V+N}_{V}$ solutions. 
Recently, the Calogero and Sutherland
quantum hamiltonians were shown to be expressible 
in terms of the operators $J^{b}_{a}$ \cite{tur3},
this result reveals the hidden symmetries of these models. 

\section{Operators preserving  P(M,V) $\oplus$ P(N,V)} 

We now  put the emphasis on the 2$\times$2 matrix 
operators which preserve the vector space
\be
P(M, V) \oplus P(N,V) \ \ \ , \ \ \ \Delta \equiv N-M \ .
\label{esp}
\ee
Without loss of generality, we assume the 
integer $\Delta$ to be non negative.
In order to classify the operators preserving (\ref{esp})
we  define a list of generators. 
First the "diagonal" generators that we choose as
\be
J^{b}_{a} (N,\Delta) =  \left(\begin{array}{cc}
J^{b}_{a}(N-\Delta) &0\\
0 &J^{b}_{a}(N)
\end{array}\right)  
- {1\over 2}   
\left(\begin{array}{cc}
1+\Delta &0\\
0        &1-\Delta 
\end{array}\right) \delta_{a}^{b} 
\label{op1}
\ee
for $0 \leq a , b \leq V$.
They are built as a direct sum of two operators 
of the type (4)-(7), translated by (\ref{trans}) in such a way that
that $J^0_0(N,\Delta)$ is proportional to the unit matrix.
The interest for this translation will appear later.

The ``non diagonal'' generators naturally split 
into ``$Q$ operators'', proportional to the matrix
$\sigma_-$  (as usual $\sigma_{\pm} = (\sigma_1 \pm i \sigma_2)/2$).
and ``$\overline Q$ operators'' proportional to the matrix
$\sigma_+$. It is convenient to write them by using a multi index 
$[A] \equiv a_1,a_2,\ldots,a_{\Delta}$.
For later convenience we also define $[\hat A_i]$
as the set $[A]$ where the index $a_i$ has been removed.
We choose the non diagonal generators respectively as
\be
Q_{[A]}
 = (-1)^{\delta} x_{a_1} \ldots x_{a_{\Delta}} \ \sigma_-  \quad ,
 0 \leq a_i \leq V  \quad , \ \ x_0 \equiv 1
\label{op2}
\ee
where $\delta$ represents the overal degree in $x_1,\dots,x_V$ 
of the monomial
$x_{a_1} \ldots x_{a_{\Delta}}$ and
\be
\overline Q^{[B]}  =
 \overline q^{[B]} \ \sigma_+ \quad ,
 0 \leq b_i \leq V  \quad 
\label{op3}
\ee
where the scalar operators $\overline q^{[B]}$, fully symmetric 
in their $\Delta$ indices  $b_k$, are defined by
\begin{eqnarray}
\overline q^{[B]} &=&\partial_{b_1} \ldots \partial_{b_{\Delta}} \quad \ \ {\rm if}
   \     0< b_1 \leq b_2 \leq b_3 \ldots \leq b_{\Delta}  \nonumber  \\
& &(D-N+\Delta-1)\partial_{b_2} \ldots \partial_{b_{\Delta}} \quad \ \ {\rm if}
   \      0= b_1 < b_2 \leq b_3 \ldots \leq b_{\Delta}   \\
& &(D-N+\Delta-1)(D-N+\Delta-2) \partial_{b_3} \ldots \partial_{b_{\Delta}} \quad \ 
    {\rm if}
     \    0= b_1 = b_2  < b_3 \ldots \leq b_{\Delta}  \nonumber \\
& &(D-N+\Delta-1)(D-N+\Delta-2) \ldots (D-N)     \quad \ 
     \ {\rm if}
      \      0= b_1 = b_2  = b_3 = b_{\Delta} \ .  \nonumber 
\end{eqnarray}
The operators $Q_{[A]}$ (and similarly the  $\overline Q^{[A]}$) 
are fully symmetric in their $\Delta$ indices $a_k$. 
Hence there are $C_{\Delta}^{V+\Delta}$ independent operators of both types. 
We then have the following proposition.

\vspace{0.5cm}
\noindent {\bf {Proposition  1}}
\vspace{0.5cm}

\noindent {\it The operators preserving the space 
$P(N-\Delta, V) \oplus P(N,V)$
are the elements of the enveloping algebra constructed over 
the generators (\ref{op1}),(\ref{op2}),(\ref{op3}) }.
\vspace{0.5cm}

This result 
(whose demonstration follows the same lines 
as in the scalar case \cite{tur1})
allows to write formally
all the operators preserving (\ref{esp}).
However, in order to classify these operators, 
it is useful to set up normal ordering rules
between the generators. 
In particular, these rules allow to write any product of operators
(the enveloping algebra) in a
canonical form, e.g. as a sum of terms where, in each term,
the $Q$ operators (if any)
are written on the left, the $\overline Q$ operators (if any) on
the right and the $J$ operators in between.
As we show next, such rules exist for the operators
(\ref{op1}),(\ref{op2}),(\ref{op3}).

 \vspace{0.5cm}
 \noindent {\bf {Normal  ordering  rules}}
 \vspace{0.5cm}

The operators (\ref{op1}) obey the commutation rules (\ref{comj})
and assemble into a reducible representation of $gl(V+1)$  
when acting on the vector space (\ref{esp}). 
The dimension of it is
\be 
C^{N+V}_V + C^{N+V-\Delta}_V  \ .
\ee

By construction, the operators $Q_{[A]}$
(resp. $\overline Q^{[A]}$)
transform as an irreducible  multiplet of dimension $C^{V+\Delta}_{\Delta}$
under the adjoint action of the generators $J^{b}_{a}(N,\Delta)$.
More precisely, we have
\bea
\label{comjq}
[ J^{b}_{a}, Q_{[A]} ] &=&
 k \  \delta^{b}_{a} Q_{[A]}
- \sum_{k=1}^{\Delta} 
\delta^{b}_{a_k} 
Q_{[\hat A_k, a]}
 \\
\label{comjqb}
[ J^{b}_{a}, \overline Q^{[A]} ] &=&
- k \ \delta^{b}_{a} \overline Q^{[A]}
+ \sum_{k=1}^{\Delta} 
\delta_{a}^{a_k} 
\overline Q^{[\hat A_k, b]} \ .
\eea
The explicit form of the generators leads to the value 
$k = \Delta$.  The first Casimir constructed with (\ref{op1}),
i.e.
\be
 T \equiv  \sum_{a=0}^V  J_a^a({N,\Delta}) \ ,
\ee
 plays the role of a grading operator :
\be
   [T,J_a^b] = 0 \ \ , \ \ 
   [T,Q_{[A]}] = \Delta V Q_{[A]} \ \ , \ \ 
   [T,\overline Q^{[A]}] = -  \Delta V \overline Q^{[A]}  \ .
\label{grading}  
\ee

The product of any two  operators $Q$ 
(and separately of two $\overline Q$'s) vanishes,
hence also their anticommutator 
\be
       \{ Q_{[A]} ,  Q_{[C]} \} = 0 \ \ \ , \ \ \  
       \{ \overline Q^{[B]} , \overline Q^{[D]} \} = 0 \ .
\label{comqq0}
\ee
The evaluation of the anti-commutator $\{ Q, \overline Q \} $ is more
involved. Its form can be guessed 
from the covariance under $gl(V+1)$, from the symmetries 
of $Q$ and $\overline Q$ in their indices and 
from the fact that the anti-commutator involves at most derivatives of
the order $\Delta$. It is therefore likely that 
the anticommutator  $\{ Q, \overline Q \} $  should be
expressed as a combination  of the tensors
\be
 W_{[A]}^{[B]}(k) \equiv {1 \over (\Delta !)^2}
               S[A] \ S[B]
    ( J_{a_1}^{b_1}  J_{a_2}^{b_2}
      \ldots  J_{a_k}^{b_k}
       \delta_{a_{k+1}}^{b_{k+1}} \ldots
       \delta_{a_{\Delta}}^{b_{\Delta}} )
\ee 
where the operator $S[.]$ denotes the sum over all permutations 
of all indices entering in the argument $[.]$.
After calculation, we found the following relations between 
(\ref{op1}),(\ref{op2}),(\ref{op3}), 
\be
    \{ Q_{[A]} , \overline Q^{[B]} \} =
    \sum_{k=0}^{\Delta} \alpha_k W_{[A]}^{[B]}(k)
\label{comqqb}
\ee
and the parameters $\alpha_k$ are numbers which 
are uniquely determined by the polynomial equation
\be
    \prod_{j=0}^{\Delta-1} (y+j) = 
    \sum_{k=0}^{\Delta} \alpha_k (y + {\Delta -1 \over 2})^k  \ .
\label{cond}
\ee
As a consequence of (\ref{cond}),
the right hand side of (\ref{comqqb}) is an even (resp. odd) polynomial
in the operators $J$ if $\Delta$ is even (resp. odd).
We would like to stress that this particularly simple expression
is due to the labelling of the generators and to the 
translation used in (\ref{op1}).
A priori, the undetermined coefficients
$\alpha_k$  could be $2 \times 2$ diagonal  matrices.

The non vanishing parameters $\alpha_k$ appear only for 
 $k={\Delta}, {\Delta-2},{\Delta-4},\ldots$ and read as follows
 for the first few values of $\Delta$~: 
\begin{eqnarray}
&\Delta = 1  \ \ , \ \ &\alpha_k : \ 1 \                         \nonumber \\ 
&\Delta = 2  \ \ , \ \ &\alpha_k : \ 1 \ , \ -{1\over 4}         \nonumber \\
&\Delta = 3  \ \ , \ \ &\alpha_k : \ 1 \ , \ -1                  \nonumber \\
&\Delta = 4  \ \ , \ \ &\alpha_k : \ 1 \ , \ -{5\over 2} \ ,\ {9\over 16} \nonumber \\
&\Delta = 5  \ \ , \ \ &\alpha_k : \ 1 \ , \ - 5 \  , \  4   \nonumber \\
&\Delta = 6  \ \ , \ \ &\alpha_k : \ 1 \ , \  -{35 \over 4} \ , \
{259 \over 16} \  , \  -{225 \over 64} \ .
\end{eqnarray}

\section{Abstract algebras}

\par We now investigate the possibility that the operators 
(\ref{op1}),(\ref{op2}),(\ref{op3}) 
represent  the generators of an abstract associative 
algebra. We will see  that there are two
types of such algebras that we note generically  ${\cal A}(V,\Delta)$
and  ${\cal B}(V,\Delta)$. 
A few  cases are known  to
coincide with Lie algebras \cite{shtu,bk}
\be
{\cal A} (1,0) \simeq sl(2) \otimes sl(2)
\ee
\be
{\cal A}(1,1) \simeq osp (2,2) \simeq spl (2,1) \ .
\ee
For $\Delta > 1$,  ${\cal A}(1,\Delta)$ corresponds
to a non linear superalgebra \cite{bk}. The algebra 
${\cal A}(1,2)$ was treated in great detail in \cite{bgkn}. 
Here we want to move away from the case $V=1$ in order to access
the hidden symmetries of the operators 
preserving (\ref{esp}) in general.

With the aim to promote the normal ordering rules 
of the previous section into a set of relations
defining an abstract associative algebra, we first note 
that the operators $J_{a}^{b}$  (resp. $Q$, $\overline Q$)  
should naturally
be interpreted as the bosonic (resp. fermionic) generators
of the algebra (this refers of course to the most natural
choice of the commutator or of the anti-commutator used to
exchange the order between these generators). 
Therefore, we expect some graded algebras to come out. 
However, it is well known that the knowledge of 
a particular representation (here
 (\ref{op1}),(\ref{op2}),(\ref{op3})) is not
sufficient in general to infer the whole algebraic structure : 
the Jacobi identities are not automatically fulfilled.
In the present case, the identities which are not obeyed
are those involving a  $\{Q,Q \}$ 
(or a  $\{\overline Q,\overline Q \}$) anticommutator
(remember that they vanish).
Although we can try to modify the whole set of (anti) commutation
relations, we limit our research of the underlying abstract
algebras in relaxing only the relation (\ref{comqq0}). 
In order to present the way to modify it, a few 
notations are worth introducing.

Due to its symmetry in the indices $[A]$, 
the representation defined by the  
 $Q_{[A]}$
(and similarly by the  $\overline Q^{[A]}$)
corresponds to a Young diagram with one line of $\Delta$ boxes.
The products
\be
  Q_{a_1 \ldots a_{\Delta} }     Q_{c_1 \ldots c_{\Delta} }
\label{rep}
\ee
assemble into a representation of  $gl(V+1)$ 
under the adjoint action of the operators $J$. 
This representation can be decomposed into irreducible pieces.
The symmetry of $Q$ is such that the irreducible representations
appearing  in the decomposition of (\ref{rep}) correspond to the 
Young diagrams consisting of two lines with total number of
$2 \Delta$ boxes.
When applied to the anticommutators 
\be
    Q_{a_1 \ldots a_{\Delta} }  Q_{c_1 \ldots c_{\Delta} }
 +  Q_{c_1 \ldots c_{\Delta} }  Q_{a_1 \ldots a_{\Delta} } \ ,
\ee
the same decomposition 
selects only the representations which are symmetric under the
exchange $[A] \leftrightarrow [C]$. In terms of Young diagrams
they correspond to the diagrams  with two lines and total number of
$2 \Delta$ boxes; the upper line is of length 
$2 \Delta-2p$ (with $2p \leq \Delta$)
 and the lower line is of even length $2p$. 
One Young tableau, corresponding to this Young diagram with fixed $p$,
is obtained by filling the first (resp. the second)  line with 
\be
[a_1,a_2,\ldots,a_{\Delta},c_{2p+1},c_{2p+2},\ldots,c_{\Delta}]
 \ \ \ , \ \ \ (\ {\rm resp.\ }
[c_1,c_2,\ldots,c_{2p}] \ ) \ .
\ee
The Young element $S_Y$ corresponding to this Young tableau reads
\be
    S_Y = 
    S[a_1, \ldots a_{\Delta},c_{2p+1}, \ldots,c_{\Delta}] 
    S[c_1, \ldots c_{2p}] E_x \ \ \ \ , \ \ \
\label{young}
\ee
with 
\be 
    E_x \equiv  \prod_{k=1}^{2p} (E - (a_k,c_k)) 
\label{exchange}
\ee
where the operator $S[.]$ (defined previously) denotes the sum
over all permutations of all indices appearing in the argument $[.]$, 
where $(a,b)$ denotes the transposition $a \leftrightarrow b$
and where $E$ is the identity operator. 

With these notations, we are ready to describe
the conditions we choose for the anticommutations of two $Q$ operators .
We restrict them by imposing
\be
    S_Y  \{ Q_{[A]} ,  Q_{[C]} \} = 0  \ .
\label{comqqy}
\ee
This corresponds to the vanishing of 
a particular representation
contained in decomposition of the symmetrized product 
of $Q_{[A]}$ with  $Q_{[C]}$ into irreducible representations
of $gl(V+1)$. 
The absent  representation is exactly the one 
related to the Young diagram defined above, 
characterized by  $p$
and by the Young element (\ref{young}). 

We have studied the associativity conditions, i.e. the
Jacobi identities, compatible with the
relation (\ref{comj}),
the relations (\ref{comjq}),(\ref{comjqb}) for an arbitrary value
of the parameter $k$, the relation 
(\ref{comqqb})  for arbitrary  values  of the  parameters
$\alpha_k$ and the relation (\ref{comqqy}) for an arbitrary integer $p$
(in fact $2p \leq \Delta$).
The results of our calculation is summarized by the following 
proposition :

\vspace{0.5cm}
\noindent {\bf {Proposition 2}}
\vspace{0.5cm}

\noindent {\it
The set of relations 
(\ref{comj}),(\ref{comjq}),(\ref{comjqb}),(\ref{comqqb}),(\ref{comqqy})
are compatible with all the  Jacobi identities
in two cases only :

1. $p=0$ and $k=\Delta$

2. $\Delta$ even, $\Delta = 2p$ and $k=-1$ }
\vspace{0.5cm}

 In both cases, the anticommutation relations
of two $\overline Q$'s have to follow the same symmetry
pattern as the anticommutation relations (\ref{comqqy})
of two $ Q$'s. 
Associativity is realized irrespectively
of the values of the parameters $\alpha_k$.
That is to say that we obtained two families of associative 
algebras ${\cal{A}},{\cal{B}}$, each indexed by $\Delta + 1$ parameters and
by the integers $V$ and $\Delta$.
By a suitable choose of the normalisation of the  $Q$'s
and/or of the $\overline Q$'s one can set 
$\alpha_{\Delta} = 1$ in (\ref{comqqb}). 
 One can also set  $\alpha_{\Delta -1}=0$ 
by using an appropriate translation (\ref{trans})
on the operators $J(N,\Delta)$.
Before presenting the proof of this result, let us
discuss a few properties of the algebras.

\vspace{0.5cm}
\begin{description}
\item [Case 1. The abstract algebra ${\cal A}(V,\Delta,\alpha_k)$] 
\ \ \ \ 

In the case $p=0$, $k = \Delta$,
the constraint (\ref{comqqy}) reads
\be
S[A,C] 
  \{ Q_{[a_1, a_2 \ldots a_{\Delta}]},
     Q_{[c_1, c_2 \ldots c_{\Delta}]}
  \}     = 0
\label{comqq1}
\ee
and the same relation has to be imposed on the operators $\overline Q$.
Using some combinatoric, one can show that (\ref{comqq1})
encodes a total number of $C^{2 \Delta + V}_V$ independent relations
among the anticommutators of two $Q$'s. Remembering that there are
$C^{V+\Delta}_V$ operators $Q$, we see easily that 
the number of constraints is lower than the number of
independent anticommutators,
that is to say that not all anticommutators are constrained.

The  operators (\ref{op1}),(\ref{op2}),(\ref{op3})
constitute a particular representation of the algebras 
of this type~: the ones corresponding to the values   
$\alpha_k$ determined by (\ref{cond}).
For these operators the conditions (\ref{comqq1}) 
are  trivially realized.

In the case $\Delta = 1$, the relation
(\ref{comqqy}) just implies that all
anticommutators of two operators $Q$ vanish 
(and similarly for  two $\overline Q$).
The algebra ${\cal A}(V,1)$ is linear, 
it coincides with the  Lie superalgebra denoted
 $spl(V+1,1)$ in the classification \cite{kac}.
If  $V=1$ one recovers the algebra $osp(2,2)$
(remember the equivalence of $osp(2,2)$ with  $spl(2,1)$).  

\vspace{0.5cm}

\item[Case 2. The abstract algebra ${\cal B}(V,\Delta,\alpha_k)$] 
\ \ \ 

If $\Delta$ is even and if $p=\Delta/2$ 
the constraints (\ref{comqqy}) on the $Q$'s can be set in the form
\be
   \lbrace Q_{a_1 a_2 \ldots a_{\Delta}},
           Q_{a_{\Delta+1} a_{\Delta +2} \ldots a_{2\Delta}}  \rbrace 
=  \lbrace Q_{\sigma(a_1) \sigma(a_2) \ldots \sigma(a_{\Delta})},
           Q_{\sigma(a_{\Delta+1}) \sigma(a_{\Delta +2})
           \ldots \sigma(a_{2\Delta})}   \rbrace 
\ee
for any permutation $\sigma$ of the $2\Delta$ indices.
The total number of independent constraints is not as easy to find
as in the case 1; we obtained it in two particular cases 
\be 
   {\Delta (\Delta -1) \over 2} \ \ \ \ {\rm if} \ \ \  V = 1
\ee
and 
\be 
   {V (V +1)(V^2+9V-4) \over 12}
 \ \ \ \ {\rm if} \ \ \  \Delta = 2 \ .
\ee
\end{description} 

\vspace{0.5cm}
\noindent {\bf {Exceptional solutions}}
\vspace{0.5cm}

It should be stressed that associative algebra could
also exist with the same structure as above, 
i.e. with (\ref{comj}), (\ref{comjq}), (\ref{comjqb}), (\ref{comqqy})
and (\ref{comqqb}) but where some of 
the parameters $\alpha_k$ are $2\times2$ diagonal matrices.
That is to say they depend on the Casimir operators
constructed with the $gl(V+1)$ subalgebra generated by the
operators (\ref{op1}).
We could not solve this problem for  general values
of $\Delta$ but we studied completely the 
cases $\Delta = 1,2,3$. We obtained one new solution 
in the case $\Delta=2,V=1,p=1$.
The  most general relation for $ \{Q, \overline Q \}$
which is compatible with associativity depends on four parameters.
It is of the form
\be
 \{ Q_{a_1 a_2} ,  \overline Q^{b_1 b_2} \} 
 = 
   \sum_{j=0}^{2} \alpha_j \ W_{a_1 a_2}^{b_1 b_2}(j)
    +
    \beta \Bigl(   C_1\  W_{a_1 a_2}^{b_1 b_2}(1)
            + (4 C_2 - 3 C_1^2) \ W_{a_1 a_2}^{b_1 b_2}(0) 
          \Bigr)
\ee
where $\beta$ is the additional  parameter 
while $C_1,C_2$ represent the Casimir operators 
(\ref{cas}) computed in the representation (\ref{op1}).     

\section{Proof of proposition 2}

Let us come to the proof of proposition 2.
The relevant Jacobi identities are 
\be
\left[\{Q_{[A]},\overline Q^{[B]}\},Q_{[C]}\right]
+\left[\{Q_{[C]},\overline Q^{[B]}\},Q_{[A]}\right]
+\left[\{Q_{[A]},Q_{[C]}\},\overline Q^{[B]}\right]=0 \ .
\ee
The application of $S_Y$ (see (\ref{young}))
to this equation and the use of (\ref{comqqy})
lead to the necessary and sufficient conditions
\be
S_Y\left (\left[\{Q_{[A]},\overline Q^{[B]}\},Q_{[C]}\right]
+\left[\{Q_{[C]},\overline Q^{[B]}\},Q_{[A]}\right]\right )=0 \ .
\ee
Moreover $S_Y$ (with an even second line) applied to a tensor 
$T_{[A,C]}$, symmetrical in $[A]$ on one side and in $[C]$ on the
other side, selects automatically the piece in $T$ symmetrical
under the exchange $[A]\leftrightarrow [C]$. Hence, the
necessary and sufficient condition becomes simply
\be
S_Y\left(\left[\{Q_{[C]},\overline Q^{[B]}\},Q_{[A]}\right]\right )=0
\ .  
\label{jacobi}
\ee

Let us first suppose that the anticommutation relations 
of $Q$ and  $\overline Q$ take the form
\be
\left\{Q_{[A]},\overline Q^{[B]}\right\}=
     S[A]S[B]J^{b_1}_{a_1}J^{b_2}_{a_2}\ldots J^{b_{\Delta}}_{a_{\Delta}}
\label{comqqbs}
\ee
rather than the more general one (\ref{comqqb}).
Using (\ref{comqqbs}) together with (\ref{comjq}),
and separating the terms, say $X'$, which come out without $k$ 
(through (\ref{comjq}))
from the terms, say  $k Y'$, which come out linear in $k$, 
the expression (\ref{jacobi}) becomes
\be
X'+ k Y'=0
\ee
where 
\bea
X'&=&-S_Y S[B]\sum_{i=1}^{\Delta}\sum_{j=1}^{\Delta}
    \delta_{a_i}^{b_1}Q_{[c_j,\hat A_i]}S[\hat C_j]
    J^{b_2}_{c_1}J^{b_3}_{c_2}\ldots J^{b_{\Delta}}_{c_{\Delta}}
\ \  - \ \ \ldots
\label{XP} \\
Y'&=& S_Y S[B]\sum_{j=1}^{\Delta}
    \delta_{c_j}^{b_1}Q_{[A]}S[\hat C_j]
    J^{b_2}_{c_1}J^{b_3}_{c_2}\ldots J^{b_{\Delta}}_{c_{\Delta}}
\ \    + \ \  \ldots  \ . 
\label{YP}
\eea
In (\ref{XP}) and (\ref{YP}), the $\dots$ refer to the terms where the
$Q$ does not appear as the first operator,
but rather after a $J$ operator. Remark
also that the index $c_j$ is absent in the set 
$[\hat C_j]$ and accordingly does not appear as a lower index in the
$J$'s.. It follows that, as it should, 
the number ($\Delta -1$) of indices $b_k$ in the product of the $J$'s
matches the number of $c_m$ indices.

Since $X'+ k Y'$ has to be zero identically, every coefficient of
every (independent) operator entering in it has to be zero. This
allows a great simplification in the necessary and sufficient
conditions. 
\vskip 0.5 cm
\begin{itemize}
\item The terms labeled $\ldots$ in (\ref{XP}) and (\ref{YP}) can
be forgotten altogether. Indeed, the terms where the $Q$'s are in
the first position are independent of the terms where they are not.
\item The symmetry on the $[B]$ can also be eliminated. Every
term, for every value of the indices $b_k$, has to vanish on its
own.
\item Let us introduce the notations 
\be
W(a_i)=\delta_{a_i}^{b_1} 
\ee
for some arbitrary fixed value of $b_1$ and
\be 
V[\hat C_k]=J^{b_2}_{c_1}J^{b_3}_{c_2}\ldots J^{b_{\Delta}}_{c_{\Delta}} 
\ee
where $c_k$ is absent as a lower index and
$b_2,\ldots,b_{\Delta}$ have also fixed values.. 
\end{itemize}
\vskip 0.5 cm 
With these simplifications, the condition $X'+kY'=0$ reduces to
the necessary and sufficient condition $X+kY=0$ with
\be
X=-S_Y \sum_{i=1}^{\Delta}\sum_{j=1}^{\Delta}
    W(a_i) Q_{[c_j,\hat A_i]}S[\hat C_j] V[\hat C_j]
\label{X}
\ee
and
\be
Y=  S_Y \sum_{j=1}^{\Delta}
    W(c_j) Q_{[A]}S[\hat C_j]
    V[\hat C_j] \ .
\label{Y}
\ee
The operator $X+k Y$ is composed of  exactly two types of 
independent operators. They can be written canonically as
\bea
O_1&=&W(c_1)Q_{[A]}V[\hat C_1]  \ ,   \\
O_2&=&W(c_{\Delta})Q_{[A]}V[\hat C_{\Delta}]  \ .
\label{OP}
\eea
Indeed~:
\vskip 0.5 cm
\begin{itemize}
\item The indices of the $Q$ operator have to be completely
symmetrical. Hence they must belong to the first line of the
Young tableau and by
symmetry of $S_Y$ can be chosen as the $[A]$ set. 
\item If the index in
$W$ is taken in the first line, it can be chosen to be $c_{\Delta}$.
This is due to the fact that any of the indices 
(except the those belonging to the set $[A]$ which already pertain
to the $Q$) in the first line is equivalent by symmetry to any
other in the first line. The remaining indices for the
$V$ can the be chosen in any order and for example in the natural order.   
\item If the index in
$W$ is taken in the second line, it can be chosen to be $c_{1}$.
Indeed any of the indices in the second line is equivalent by symmetry to any
other in the second line. The remaining indices for the
$V$ can again be chosen in any order and for example in the natural order.   
\end{itemize}
\vskip 0.5 cm
The remaining task is to extract in $X$ and in $Y$ the number of times
the operators $O_1$ and $O_2$ occur. This is a rather delicate
operation in terms of the symmetries involved.
Let us call $X_i$ (resp.  $Y_i$) with $i=1,2$ the coefficient of
the operator $O_i$ in $X$ (resp. $Y$) 
With these notations 
the condition $X+kY=0$ becomes equivalent to 
\be
    X_1 + k Y_1 = 0 \ \ \ , \ \ \ X_2 + k Y_2 = 0 \ .
\label{xi+yi}
\ee

To now compute these four coefficients, we will make
use of the 
fundamental theorem of group theory which states that,  
if $P$ is any permutation of the elements in $[A]$ 
\be
S[A]=PS[A]=S[A]P   
\label{fonda}
\ee
\vskip 0.5 cm
\begin{description}
\item[Computation of $X_1$] 
\ \ \ 

Let us rewrite $X$ as
\bea
X=&-&
   \sum_{i=1}^{\Delta}\sum_{j=1}^{\Delta}
   S[A,c_{2p+1},\ldots,c_{\Delta}] S[c_1,\ldots,c_{2p}]
   E_x   \nonumber \\
   && W(a_i) Q_{[c_j,\hat A_i]}S[\hat C_j] V[\hat C_j]
\label{X1a}
\eea
where we have interchanged the finite sommation on $i$ and $j$
with the symmetry operations.
First, the $a_i$ in $W(a_i)$ which belongs to the first line has
to be replaced by a $c$ belonging to the second line. This can
be done at the intervention of the operator $E_x$ only. At the
same time none of the other $a_j,j \neq i$ in $Q$ should be
replaced by an element of the second line. Hence, from the
$2^{2p}$ terms in $E_x$ we can restrict ourselves
to the transposition $(a_i,c_i)$ which comes
with a minus sign. At the same time $i$ can be restricted to the range $1,2p$.
The sommation on $j$ then has one term with $j=i$. For the terms
with $j \neq i$, the $c_j$ in $Q$ has to belong to the set
$j=2p+1,\ldots,\Delta$ in order to be able to replace it by an
$a$ by the first symmetry operator $S$ in (\ref{X1a}). Hence
the restricted part of $X$, say $\hat X$, is composed of two
pieces, say $\hat X_{\alpha}$ and $\hat X_{\beta}$, 
\bea
\hat X_{\alpha}=&&
   \sum_{i=1}^{2p}
   S[A,c_{2p+1},\ldots,c_{\Delta}] S[c_1,\ldots,c_{2p}]
   (a_i,c_i)
   W(a_i) Q_{[c_i,\hat A_i]}S[\hat C_i] V[\hat C_i]   
         \nonumber  \\
          =&&
   \sum_{i=1}^{2p}
   S[A,c_{2p+1},\ldots,c_{\Delta}] S[c_1,\ldots,c_{2p}]
   W(c_i) Q_{[A]}S[\hat C_i] V[\hat C_i]              
         \nonumber  \\
          =&&
   \sum_{i=1}^{2p}
   S[A,c_{2p+1},\ldots,c_{\Delta}] S[c_1,\ldots,c_{2p}]
   (c1,c_i)
    W(c_i) Q_{[A]}S[\hat C_i] V[\hat C_i]         
         \nonumber  \\
          =&&
   S[A,c_{2p+1},\ldots,c_{\Delta}] S[c_1,\ldots,c_{2p}]
   \sum_{i=1}^{2p}
     W(c_1) Q_{[A]}S[\hat C_1] V[\hat C_1] 
\label{X1alpha}
\eea
and
\bea
\hat X_{\beta}=&&
   \sum_{i=1}^{2p}\sum_{j=2p+1}^{\Delta}
   S[A,c_{2p+1},\ldots,c_{\Delta}] S[c_1,\ldots,c_{2p}]
   (a_i,c_i)
    \nonumber  \\
   && W(a_i) Q_{[c_j,\hat A_i]}S[\hat C_j] V[\hat C_j]  \ .
\label{X1beta}
\eea

Using (\ref{fonda}), we
easily conclude that in $X_{\alpha}$ the following coefficient appears 
\be
(\Delta)!(2p)!(\Delta-2p)!  \ .
\label{xbeta}
\ee
The first factor $(\Delta)!$ comes from the permutation of the $[A]$ set
which always contributes to an equal factor due to the symmetry
of the $Q$. The second term $(2p)!$ comes from the product of
 of the sommation over $i$ (a factor $2p$)
and of a factor $(2p-1)!$ coming from the repetition
of the symmetries in $[c_2,\ldots,c_{2p}]$ contained in the
second and in the third $S$ factors. The last term $(\Delta-2p)!$ 
comes from from the repetition
of the symmetries in $[c_{2p+1},\ldots,c_{\Delta}]$ contained in the
first and in the third $S$ factors.

Let us now focuss our attention on the $X_{\beta}$ term. Using
(\ref{fonda}) we can factor out of
$S[A,c_{2p+1},\ldots,c_{\Delta}]$, at no cost, a transposition
factor $(a_i,c_j)$ and from $S[c_1,\ldots,c_{2p}]$ a factor
$(c_1,c_i)$. The product of these two transpositions together
with the transposition in $X_{\beta}$ leads to the cyclic
permutation $(a_i,c_1,c_i,c_j)$ and the relevant part $\hat X_{\beta}$
becomes
\bea
\hat X_{\beta}=&&
   \sum_{i=1}^{2p}\sum_{j=2p+1}^{\Delta}
   S[A,c_{2p+1},\ldots,c_{\Delta}] S[c_1,\ldots,c_{2p}]
     \nonumber  \\
    &&W(c_1) Q_{[A]}S[\hat C_1] V[\hat C_1] \ .
\label{X1beta2}
\eea

The following coefficient then appears
\be
(\Delta-2p)(\Delta)!(2p)!(\Delta-2p)!  \ .
\label{xalpha}
\ee
The extra factor as compared to the coefficient coming out of 
$X_{\alpha}$ is due to the extra sommation over $j$.

Summing up the results (\ref{xalpha},\ref{xbeta}), we find
\be
X_1=(1+\Delta-2p)(\Delta)!(2p)!(\Delta-2p)! \ \ .  
\label{X_1}
\ee
\vskip 0.5 cm
\item[Computation of $Y_1$]
\ \ \ 

The same technique applied to $Y_1$ is much simpler as the
relevant term in $E_x$ is simply the identity. Hence
\bea
\hat Y=&&
     S[A,c_{2p+1},\ldots,c_{\Delta}] S[c_1,\ldots,c_{2p}]
     \sum_{j=1}^{2p}
     W(c_j) Q_{[A]}S[\hat C_j]
    V[\hat C_j]   \nonumber  \\
=&&
     S[A,c_{2p+1},\ldots,c_{\Delta}] S[c_1,\ldots,c_{2p}]
     \sum_{j=1}^{2p}
     W(c_1) Q_{[A]}S[\hat C_1]
    V[\hat C_1]  \ .
\label{Y1}
\eea
To pass from the first to the second line we have factored out of
$S[c_1,\ldots,c_{2p}]$ the transposition $(c_1,c_j)$.

Collecting again the factors, we find
\be
Y_1=(\Delta)!(2p)!(\Delta-2p)! \ \ .  
\label{Y_1}
\ee
\vskip 0.5 cm
\item[Computation of $X_2$]
\ \ \ 

The relevant term in $E_x$ is again the identity and
the reduced part of $X$ which can lead to a term of the form $O_2$
(\ref{OP}) is
\be
\hat X=-   \sum_{i=1}^{\Delta}\sum_{j=1}^{\Delta}
     S[A,c_{2p+1},\ldots,c_{\Delta}] S[c_1,\ldots,c_{2p}]
    W(a_i) Q_{[c_j,\hat A_i]}S[\hat C_j] V[\hat C_j] \ .
\label{X2a}
\ee
The sommation on $j$ has to be restricted to those values in the
first line of the Young diagram. A transposition $(a_i,c_j)$ can
then be factored out of $S[A,c_{2p+1},\ldots,c_{\Delta}]$ as
well as a
transposition $(a_i,c_{\Delta})$, i.e. in total a cyclic permutation
$(a_i,c_{\Delta},c_j)$ . We find
\bea
\hat X=&-& \sum_{i=1}^{\Delta}\sum_{j=2p+1}^{\Delta}
     S[A,c_{2p+1},\ldots,c_{\Delta}] S[c_1,\ldots,c_{2p}]
     (a_i,c_{\Delta},c_j)
    \nonumber  \\
    && W(a_i) Q_{[c_j,\hat A_i]}S[\hat C_j] V[\hat C_j]   
     \nonumber  \\
     =&-& \sum_{i=1}^{\Delta}     
     S[A,c_{2p+1},\ldots,c_{\Delta}] S[c_1,\ldots,c_{2p}]
     \sum_{j=2p+1}^{\Delta}
    \nonumber  \\
    && W(c_{\Delta}) Q_{[A]}S[\hat C_{\Delta}] V[\hat C_{\Delta}] \ .
\label{X2b}
\eea
Collecting the factors as usual, we find
\be
X_2=-\Delta(\Delta)!(2p)!(\Delta-2p)! \ .
\label{X_2}
\ee
The extra factor ($\Delta$) comes from the sommation over $i$.
\vskip 0.5 cm
\item[Computation of $Y_2$]
\ \ \ 

In this last case the relevant term in $E_x$ is again
the identity and
the reduced part of $Y$ which can lead to a term of the form $O_2$
(\ref{OP}) is
\bea
\hat Y&=&
     S[A,c_{2p+1},\ldots,c_{\Delta}] S[c_1,\ldots,c_{2p}]
     \sum_{j=2p+1}^{\Delta}
    W(c_j) Q_{[A]}S[\hat C_j]
    V[\hat C_j]  \nonumber        \\
    &=&
     S[A,c_{2p+1},\ldots,c_{\Delta}] S[c_1,\ldots,c_{2p}]
     \sum_{j=2p+1}^{\Delta} (c_j,c_{\Delta})
    W(c_j) Q_{[A]}S[\hat C_j]
    V[\hat C_j]    \nonumber   \\
    &=&
     S[A,c_{2p+1},\ldots,c_{\Delta}] S[c_1,\ldots,c_{2p}]
     \sum_{j=2p+1}^{\Delta}
    W(c_{\Delta}) Q_{[A]}S[\hat C_{\Delta}]
    V[\hat C_{\Delta}]        \ .
\label{Y2}
\eea
Collecting the terms, we find
\be
Y_2=(\Delta)!(2p)!(\Delta-2p)!  \ .
\label{Y_2}
\ee
\end{description}
\vskip 0.5 cm
We can now summarize the conditions coming from (\ref{xi+yi}) :
\begin{enumerate}
\item The conditions coming from the Jacobi identities are thus two in
number if $p \neq 0$ (the condition for the operator $O_1$ to be
defined) and if $\Delta \neq 2p$ (the condition for $0_2$ to be
defined). These conditions
\bea
k&=&-(\Delta+1-2p)  \\
k&=&\Delta
\label{gencase}
\eea
are incompatible.
\item More generally, the anticommutators of two $Q$'s cannot
vanish for more than one representation.
\item If $p=0$ the only condition comes from the $0_2$ operator.
It is
\be
k=\Delta
\ee
which is a solution to our problem. The corresponding Young
diagram has only one line of length $2\Delta$.
\item If $\Delta$ is even and $\Delta=2p$ the only condition
comes from the $0_1$ operator.
It is
\be
k=-1
\ee
which is a second solution to our problem.
The corresponding Young diagram has two lines of equal length $\Delta$.
\end{enumerate}
\vskip 0.5 cm
This achieves the proof of the proposition when the 
anticommutator $\{ Q , \overline Q \}$ is restricted to 
(\ref{comqqbs})
It is easy to see that the other allowed terms in the
anticommutator of the $Q$'s with the $\overline Q$'s, i.e. those
which do not involve $J$'s only but the products of $J$'s and
$\delta$'s as in (\ref{comqqb}) lead to exactly the same restrictions.
Hence they can all be present at the same time 
leaving us with the form (\ref{comqqb}) with the  $\Delta+1$
arbitrary coefficients.

The conditions coming from the Jacobi identities involving two
$\overline Q$ and one $Q$ also lead to exactly the same conditions.
Hence the anticommutators which are chosen to be zero for the
anticommutations
of two $Q$'s on one side or of two $\overline Q$'s on the other
side must be identical.

\section{Summary and Conclusions}

The operators preserving globally a system of two polynomials in
$V$ variables ($V\geq 1$) and of degrees $N$ and $N-\Delta$
($\Delta\geq 0$) respectively can
be constructed as the elements of the enveloping algebra of certain
superalgebras. 

In this paper, we have constructed a family
of such associative, non linear superalgebras. Any of these
algebras is specified by $V$, by $\Delta$ and by a set of
$\Delta+1$ complex numbers noted $\alpha_k$ with
$k=0,1,\ldots,\Delta$. They are generated by $(V+1)^2$ (bosonic) operators
\be
J_a^b\ \ \ , \ \ \ a,b=0,1,\ldots,\Delta
\ee
and by two sets of $C_{\Delta}^{V+\Delta}$ (fermionic) operators
\be
Q_{[a_1,\ldots,a_{\Delta}]}\ \ \ , \ \ \ 
    {\overline Q}^{[a_1,\ldots,a_{\Delta}]} \ \ \ ,\ \ \ 
        a_k=0,1,\ldots,\Delta
\ee
symmetric in their $\Delta$ indices.

The bosonic generators obey the commutation relations of the Lie
algebra $gl(V+1)$. The operators $Q$ (and separately the
$\overline Q$) assemble into a specific representation of
$gl(V+1)$ under the adjoint action of $J_a^b$ 
(see (\ref{comjq}),(\ref{comjqb})).
The anticommutators $\{Q,\overline Q\}$ are polynomials of
degree at most $\Delta$ in the bosonic operators. The
arbitrariness of the polynomials is coded in the 
$\Delta+1$~parameters $\alpha_k$ (\ref{comqqb}).

All the supplementary conditions on the products of the
operators $Q$ (and of the operators $\overline Q$) necessary to
guarantee associativity (equivalent to the generalised Jacobi
identities) are given by our proposition~2.

For all fixed values of the integers $V$ and $\Delta$ and of the
complex parameters $\alpha_k$ we denote
${\cal{A}}(V,\Delta,\alpha_k)$ the algebra corresponding to 
case~1 of proposition~2. If $\Delta$ is even, a supplementary
algebraic structure, that we denote
${\cal{B}}(V,\Delta,\alpha_k)$ is also possible, as predicted by
case~2 of proposition~2.

Refering to the general definition of a $W$-algebra given
recently in \cite{sorba}, it is natural to classify 
${\cal{A}}(V,\Delta,\alpha_k)$ and ${\cal{B}}(V,\Delta,\alpha_k)$ 
as ``finite $W_{\Delta+1}$-superalgebras''.

An analysis of the representations of
${\cal{A}}(1,2,\alpha_0, \alpha_1, \alpha_2)$, performed
recently \cite{bgkn}, leads to a rather rich set of inequivalent
irreducible, finite dimensional representations .

Let us stress again that the operators in the enveloping
algebras that we have constructed are directly relevant for the
study of quasi-exactly solvable systems of equations.

\vfill\eject

\end{document}